\documentclass[a4paper,11pt,reqno]{amsart}
\usepackage{amsxtra}
\usepackage{color}
\usepackage{amsopn}
\usepackage{amsmath,amsthm,amssymb,arydshln}
\usepackage{mathrsfs,mathtools}
\usepackage{stmaryrd}
\usepackage{graphicx}
\usepackage{bm}
\usepackage{enumerate}
\usepackage{xcolor}
\usepackage{pifont}
\usepackage{adjustbox}
\usepackage{tikz}
\usepackage{color}
\usepackage{hyperref}

\newtheorem{theorem}{Theorem}[section]

\newtheorem*{theorem*}{Theorem}

\theoremstyle{definition}

\newtheorem{example}[theorem]{Example}

\newcommand{\R}{{\mathbb R}}
\renewcommand{\H}{\mathrm{H}}

\newcommand{\beq}{\begin{equation}}
\newcommand{\eeq}{\end{equation}}

\newcommand{\f}{\varphi}

\newcommand{\G}{{\mathrm G}}

\newcommand{\W}{\wedge}

\DeclareMathOperator\tr{tr}

\DeclareMathOperator\End{End}
\DeclareMathOperator\Aut{Aut}

\DeclareMathOperator{\Der}{Der}

\newcommand{\Ric}{{\rm Ric}}

\newcommand{\frg}{\mathfrak{g}}

\newcommand{\frh}{\mathfrak{h}}
\newcommand{\frk}{\mathfrak{k}}

\newcommand{\frn}{\mathfrak{n}}

\newcommand{\frr}{\mathfrak{r}}

\newcommand{\gsl}{\mathfrak{sl}}

\newcommand{\diag}{{\rm diag}}

\newcommand{\sst}{\scriptscriptstyle}

\numberwithin{equation}{section}

\title[Remarks on homogeneous solitons of the G$_2$-Laplacian flow]{Remarks on homogeneous solitons of the G$_{\mathbf2}$-Laplacian flow}

\author{Anna Fino} 
\address{Dipartimento di Matematica ``G.~Peano'' \\ Universit\`a degli Studi di Torino\\
Via Carlo Alberto 10\\
10123 Torino\\ Italy}
\email{annamaria.fino@unito.it}
\email{alberto.raffero@unito.it}
\author{Alberto Raffero}

\begin{document}
\begin{abstract} 
We show the existence of expanding solitons of the G$_2$-Laplacian flow on non-solvable Lie groups, and we give the  
first example of a steady soliton that is not an extremally Ricci pinched G$_2$-structure. 
\end{abstract}
\maketitle

\section{Introduction}
A G$_2$-structure on a 7-manifold $M$ is said to be {\em closed} if the defining positive 3-form $\f\in\Omega^3_{\sst+}(M)$ satisfies the equation $d\f=0$. 
The intrinsic torsion of a closed G$_2$-structure can be identified with the unique 2-form $\tau$ for which $d*_\f\f = \tau\W\f = -*_\f\tau,$ 
where $*_\f$ denotes the Hodge operator associated with the Riemannian metric $g_\f$ and orientation $dV_\f$ induced by $\f$ (cf.~\cite{Bry}). 
When this {\em intrinsic torsion form} $\tau$ vanishes identically, the G$_2$-structure is called {\em torsion-free}. 

A closed G$_2$-structure $\f$ is called a {\em Laplacian soliton} if it satisfies the equation
\begin{equation}\label{LS}
\Delta_\f\f=\lambda\f +\mathcal{L}_{\sst X}\f
\end{equation}
for some real number $\lambda$ and some vector field $X$ on $M,$ where $\Delta_\f$ denotes the Hodge Laplacian of the metric $g_\f$.  
It is known that a closed G$_2$-structure $\f$ satisfies \eqref{LS} if and only if the solution of the {\em Laplacian flow} 
\begin{equation}\label{LF}
\begin{dcases}
\frac{\partial}{\partial t} \f(t) = \Delta_{\f(t)}\f(t),\\
d\f(t)=0, 
\end{dcases}
\end{equation}
starting from it at $t=0$ is self-similar (cf.~\cite{Lin,LoWe1}, and see the recent survey \cite{Lot} for more information on the Laplacian flow \eqref{LF}). 

Depending on the sign of $\lambda$, a Laplacian soliton is called {\em shrinking} ($\lambda<0$), {\em steady} ($\lambda=0$), or {\em expanding} ($\lambda>0$). 
By \cite{Lin,LoWe1}, on a compact manifold every Laplacian soliton which is not torsion-free must satisfy \eqref{LS} with $\lambda>0$ and $\mathcal{L}_{\sst X}\f\neq0.$ 
The existence of non-trivial Laplacian solitons on compact manifolds is still an open problem. 

The non-compact setting is less restrictive, and examples of Laplacian solitons of any type occur. 
Moreover, all examples obtained so far are given by seven-dimensional, connected, simply connected solvable Lie groups $\G$ endowed with a left-invariant closed $\G_2$-structure $\f$ satisfying \eqref{LS} 
with respect to a special vector field $X$ \cite{FFR,FiRa1,Lau1,Lau2,LaNi,Nic}. 
In detail, the vector field $X$ is defined by the one-parameter group of automorphisms $F_t\in\Aut(\G)$ such that $\left.dF_t\right|_{1_{\G}} = \exp(tD),$ where $D$ is a derivation of the Lie algebra $\frg$ of $\G$.  
According to \cite{Lau1}, these Laplacian solitons are said to be {\em semi-algebraic}, and they are called {\em algebraic} if the adjoint of $D$ with respect to the 
inner product $g_\f$ on $\frg$ is also a derivation. 
Notice that, even if the Lie group $\G$ admits a co-compact discrete subgroup $\Gamma$, a (semi)-algebraic soliton does not define a Laplacian soliton on the compact quotient $\Gamma\backslash\G$, 
as the vector field defined by $F_t$ does not descend to $\Gamma\backslash\G$. 

The purpose of this note is to discuss new homogeneous examples of Laplacian solitons that differ in some aspects from those appearing in the literature. 

In Section \ref{SecNS}, we show that algebraic Laplacian solitons also exist on non-solvable Lie groups. 
This gives a further evidence of the difference between homogeneous Laplacian solitons and Ricci solitons, in addition to the results obtained in \cite{Lau1, Nic}.  
Indeed, every homogeneous Ricci soliton is isometric to an algebraic soliton of the Ricci flow \cite{Jab}, 
and all known examples are isometric to a left-invariant algebraic Ricci soliton on a simply connected solvable Lie group (see \cite{Jab,LaLa} for more details).

In Section \ref{SecSLS}, we focus on steady Laplacian solitons. Currently, all known examples \cite{Lau2,LaNi} are given by 
{\em extremally Ricci pinched} $\G_2$-structures, namely closed $\G_2$-structures $\f$ whose intrinsic torsion form $\tau$ satisfies the equation $d\tau = \frac16\, |\tau|^2\,\f +\frac16\,*_\f(\tau\W\tau)$. 
In fact, by \cite{LaNi} every left-invariant extremally Ricci pinched $\G_2$-structure on a (necessarily solvable) Lie group is a steady algebraic soliton. 
Here, we show that the converse of this result does not hold. 
In detail, we give an example of a simply connected solvable Lie group endowed with a steady algebraic soliton that is not an extremally Ricci pinched $\G_2$-structure. 
This example satisfies a further remarkable property: there exists a left-invariant vector field $X$ on the Lie group for which $\Delta_\f\f = \mathcal{L}_{\sst X}\f$. 
To our knowledge, this is the first example of a left-invariant closed $\G_2$-structure satisfying the equation \eqref{LS} with respect to a left-invariant vector field.

\section{Algebraic Laplacian solitons on Lie groups}
In this section, we briefly recall some known facts on left-invariant $\G_2$-structures on Lie groups and on algebraic solitons of the $\G_2$-Laplacian flow. 
We refer the reader to \cite{FiRa,Lau1} and the references therein for more details. 

Let $\G$ be a seven-dimensional, connected, simply connected Lie group.  
It is well-known that there is a one-to-one correspondence between left-invariant $\G_2$-structures on $\G$ and $\G_2$-structures on the corresponding Lie algebra $\frg$. 

Recall that a 3-form $\f\in\Lambda^3(\frg^*)$ defines a G$_2$-structure on $\frg$ if and only if the symmetric bilinear map
\[
b_\f : \frg \times \frg \rightarrow \Lambda^7(\frg^*),\quad b_\f(u,v) \coloneqq \frac16\,\iota_u\f\W\iota_v\f\W\f,
\]
satisfies the following conditions
\begin{enumerate}[i)]
\item $\det(b_\f)^{1/9}\in\Lambda^7(\frg^*)$ is not zero; 
\item the symmetric bilinear form $\det(b_\f)^{-1/9}\, b_\f : \frg\times\frg\rightarrow\R$ is positive definite.
\end{enumerate}
When this happens, the inner product and orientation induced by $\f$ on $\frg$ are given by $g_\f = \det(b_\f)^{-1/9}\, b_\f$ and $dV_\f=\det(b_\f)^{1/9}$, respectively. 

A $\G_2$-structure $\f$ on $\frg$ is {\em closed} if it satisfies the equation $d\f=0$, where $d$ denotes the Chevalley-Eilenberg differential of $\frg$. 
Using left multiplication, it is possible to extend a closed $\G_2$-structure on $\frg\cong T_{1_\G}\G$ to a left-invariant closed $\G_2$-structure on $\G$. 

Let $\f$ be a closed $\G_2$-structure on $\frg$, and let $\tau\in\Lambda^2(\frg^*)$ be the corresponding intrinsic torsion form. 
By \cite{Lau1}, the following conditions are equivalent 
\begin{enumerate}[1)]
\item there exist a real number $\lambda$ and a derivation $D\in\Der(\frg)$ for which 
\begin{equation}\label{AS}
\Ric(g_\f) -\frac{1}{12}\, \tr(\tau_\f^2)\,\mathrm{Id} + \frac{1}{2}\,\tau_\f^2 = -\frac{\lambda}{3}\,\mathrm{Id} - D, 
\end{equation}
where $\Ric(g_\f)\in\End(\frg)$ is the Ricci endomorphism of $g_\f$, and $\tau_\f\in\End(\frg)$ is the skew-symmetric endomorphism defined via the identity $\tau(\cdot,\cdot)=g_\f(\tau_\f\cdot,\cdot)$;
\item the left-invariant closed $\G_2$-structure induced by $\f$ on $\G$ satisfies the equation
\[
\Delta_\f\f = \lambda\f + \mathcal{L}_{\sst X_{\sst D}}\f,
\]
where $X_{\sst D}$ is the vector field on $\G$ defined by the unique one-parameter group of automorphisms $F_t\in\Aut(\G)$ such that $\left.dF_t\right|_{1_{\G}} = \exp(tD)\in\Aut(\frg)$. 
\end{enumerate}
A closed $\G_2$-structure satisfying any of the above conditions is called an {\em algebraic soliton},  and it clearly defines a Laplacian soliton on $\G$.

\section{Laplacian solitons on non-solvable Lie groups}\label{SecNS}
In \cite{FiRa}, we obtained the classification of seven-dimensional, unimodular, non-solvable Lie algebras admitting closed G$_2$-structures, 
showing that only four Lie algebras of this type exist, up to isomorphism. 
In this section, we give an example of an expanding algebraic soliton on two of them. Before doing this, we explain some conventions that we will use.    

Given a Lie algebra $\frg$ of dimension $n$, we write its structure equations with respect to a basis of $1$-forms $\{e^1,\ldots,e^n\}$  by 
specifying the $n$-tuple $(de^1,\ldots,de^n)$. Moreover, we use the shorthand $e^{ijk\cdots}$ to denote the wedge product  $e^i\W e^j\W e^k\W\cdots$. 
Finally, we denote by $(e_1,\ldots,e_n)$ the basis of $\frg$ with dual basis $(e^1,\ldots,e^n)$, and we write the matrix associated with any endomorphism of $\frg$ with respect to this basis.

\begin{example}
Consider the one-parameter family of pairwise non-isomorphic, unimodular, non-solvable Lie algebras $\frg_\mu$ with the following structure equations 
\[
\frg_\mu = \left(e^{45}, -2\mu\,e^{27}, -2\,e^{37}, 2\,e^{15}, e^{14}, 2\left(1+\mu\right) e^{67}, 0 \right),\quad -1<\mu\leq-\frac12. 
\]
For each $\mu\in\left(-1,-\frac12\right]$, $\frg_\mu$ is isomorphic to the decomposable Lie algebra $\gsl(2,\R)\oplus\frr_{4,\mu,-1-\mu}$ appearing in \cite[Main Theorem]{FiRa}.

The following 3-form defines a closed G$_2$-structure on $\frg_\mu$
\[
\f_\mu  = -2\mu\,e^{127} + e^{347} + 2\left(1+\mu\right) e^{567} + e^{135} - e^{146} - e^{236} - e^{245}. 
\]
The inner product and the volume form induced by $\f_\mu$ on $\frg_\mu$ are given by 
\begin{eqnarray*}
g_{\f_\mu}		&=& \left(\frac{2\mu^2}{1+\mu}\right)^{1/3} \left[(e^1)^2+(e^2)^2\right] + \left(\frac{1}{-4\mu(1+\mu)}\right)^{1/3} \left[(e^3)^2+(e^4)^2\right] \\
			& & \left(\frac{2(1+\mu)}{-\mu}\right)^{1/3}  \left[(e^5)^2+(e^6)^2\right]  +\left(4\mu(1+\mu)\right)^{2/3} (e^7)^2, \\
dV_{\f_\mu}	&=& \left(-4\mu(1+\mu)\right)^{1/3} e^{1234567}, 	
\end{eqnarray*}
and the intrinsic torsion form of $\f_\mu$ is
\[
\tau_\mu = - 2 \left(\frac{2\mu^2}{1+\mu}\right)^{1/3} e^{12}  + \left(\frac{2}{\mu+\mu^2}\right)^{2/3}\, e^{34}   - 2 \left(\frac{2(1+\mu)^2}{\mu}\right)^{2/3} e^{56}. 
\]


Using these data, it is possible to check that the equation \eqref{AS} is satisfied for the following value of $\lambda$
\[
\lambda = 2\,(\mu^2+\mu+1)  \left(\frac{2}{\mu(1+\mu)} \right)^{2/3}, 
\]
and the following derivation of $\frg_\mu$
\[
D = -2\, \diag\left(0, \left(\frac{4(1+\mu)}{\mu^2}\right)^{1/3},   (4\mu(1+\mu))^{1/3}, 0, 0,  -\left(\frac{4\mu}{(1+\mu)^2}\right)^{1/3}, 0 \right).
\]
Since $-1<\mu\leq-\frac12$, we have that $\lambda>0$.  
Thus, $\f_\mu$ is an expanding algebraic soliton, and it gives rise to a left-invariant expanding Laplacian soliton on the simply connected non-solvable Lie group with Lie algebra $\frg_\mu$.

\end{example}

\begin{example}
Consider the one-parameter family of pairwise non-isomorphic, unimodular, non-solvable Lie algebras 
\[
\frk_\alpha = \left(-\alpha\,e^{17}, -2\,e^{35}, -e^{25}-e^{57}, \frac{\alpha}{2}\, e^{47} - e^{67}, e^{23}+e^{37}, e^{47} + \frac{\alpha}{2}\, e^{67}, 0 \right),\quad \alpha>0.
\]
This family of Lie algebras is isomorphic to the family $\gsl(2,\R)\oplus\frr'_{4,\alpha,-\alpha/2}$ appearing in \cite[Main Theorem]{FiRa}. 

The 3-form 
\[
\f_\alpha =  \frac{\alpha}{2} \left(e^{127} + e^{347} + e^{567}\right) + e^{135} - e^{146} - e^{236} - e^{245}
\]
defines a closed G$_2$-structure on $\frk_\alpha$ inducing the inner product 
\[
g_{\f_\alpha} = \sum_{i=1}^6 (e^i)^2 + \frac{\alpha^2}{4}\,(e^7)^2,
\]
and the volume form $dV_{\f_\alpha} = \frac{\alpha}{2}\,e^{1234567}.$ The intrinsic torsion form of $\f_\alpha$ is 
\[
\tau_\alpha = 4\,e^{12} - 2\,e^{34} -2\, e^{56}. 
\]

It is now possible to check that the equation \eqref{AS} is satisfied for $\lambda=6$ and 
\[
D = \diag\left(-2,0,0,4,0,4,0 \right) \in \Der(\frk_\alpha). 
\]
Thus, the 3-form $\f_\alpha$ is an algebraic soliton on $\frk_\alpha$, and it induces a left-invariant expanding Laplacian soliton on the simply connected non-solvable Lie group with Lie algebra $\frk_\alpha$. 
\end{example}

\section{A steady soliton that is not extremally Ricci pinched}\label{SecSLS}
A closed G$_2$-structure $\f$ whose intrinsic torsion form $\tau$ satisfies the equation 
\[
d\tau = \frac16\, |\tau|^2\,\f +\frac16\,*_\f(\tau\W\tau)
\]
is called {\em extremally Ricci pinched} ({\em ERP} for short). 
This class of G$_2$-structures was introduced by Bryant in \cite{Bry}. 
Homogeneous examples are discussed in \cite{Bry,FiRa2,Lau2,LaNi}, and the existence of non-homogeneous examples has also been established \cite{Bal}. 
In \cite{FiRa2}, we proved that the solution of the Laplacian flow \eqref{LF} starting from an ERP G$_2$-structure exists for all times and stays ERP. 
However, on compact manifolds ERP G$_2$-structures  cannot be steady solitons. This is not true anymore in the non-compact setting. 
Indeed, in the recent work \cite{LaNi} the authors proved that any left-invariant ERP G$_2$-structure on a simply connected Lie group is always a steady soliton. 
In this section, we show that the converse of this result does not hold.  

Let us consider the solvable Lie algebra $\frh$ with the following structure equations
\[
\frh = \left(0, 0, -e^{37}, e^{47}, 2\,e^{14}+e^{57}, -2\,e^{24} + e^{67}, 0 \right). 
\]
This Lie algebra is not unimodular and it is isomorphic to the semidirect product $\R \ltimes \frn,$ 
where $\R=\langle e_7\rangle$ and $\frn=\langle e_1,\ldots,e_6\rangle$ is a six-dimensional decomposable nilpotent Lie algebra. 

The 3-form 
\[
\f = e^{127} + e^{347} + e^{567} + e^{135} - e^{146} - e^{236} - e^{245}
\]
defines a closed G$_2$-structure on $\frh$ inducing the inner product $g_\f = \sum_{i=1}^7(e^i)^2$ and the volume form $dV_\f=e^{1234567}$. 
The intrinsic torsion form of $\f$ is
\[
\tau = 2\,e^{12} + 2\,e^{34} -4\,e^{56}, 
\]
and it satisfies 
\[
d\tau = \Delta_\f\f = -8\left(e^{146}+e^{245}-e^{567}\right).
\] 
We immediately see that $\f$ is not ERP, 
since $*_\f(\tau\W\tau) = -16\left(e^{127}+e^{347}\right) +8\,e^{567}$ and  $|\tau|^2 = 24$. 
 

The closed G$_2$-structure $\f$ is a steady algebraic soliton,  as it satisfies the equation \eqref{AS} with $\lambda=0$ and 
\[
D = \diag\left(0,0,-4,4,4,4,0\right)\in\Der(\frh). 
\]
In particular, the left-invariant closed G$_2$-structure induced by $\f$ on the simply connected solvable Lie group $\H$ with Lie algebra $\frh$ satisfies 
\[
\Delta_\f\f = \mathcal{L}_{X_D}\f.
\]

This example has a further remarkable property, since there exists a left-invariant vector field $X$ on $\H$ for which the left-invariant steady soliton $\f$ satisfies the equation
\[
\Delta_\f\f = \mathcal{L}_{\sst X}\f = d(\iota_{\sst X}\f). 
\]
In detail, $X$ is the left-invariant vector field on $\H$ induced by the vector $-4\,e_7\in\frh$.

\section*{Acknowledgements}
The authors were supported by the ``National Group for Algebraic and Geometric Structures, and their Applications'' (GNSAGA -- INdAM), 
and by the project PRIN 2015 ``Real and Complex Manifolds: Geometry, Topology and Harmonic Analysis".

\end{document}